\documentclass[a4paper]{article}
\usepackage{epsf,epsfig}
\usepackage{graphics,color}
\usepackage{amsmath}
\usepackage{amssymb}
\usepackage{cite}
\usepackage{verbatim}
\usepackage{float}
\usepackage{graphicx}
\usepackage{amsthm}
\usepackage{textcomp}
\usepackage{subfig}

\newcommand{\no}{\nonumber}

\newcommand{\calV}{{\mathcal V}}

\newcommand{\calJ}{{\mathcal J}}

\newcommand{\duavw}[1]{\langle\!\langle{#1}\rangle\!\rangle}

\newcommand{\R}{\mathbb{R}}
\newcommand{\RR}{\mathbb{R}}

\newcommand{\Om}{\Omega}

\def\eps{\epsilon}

\def\Om{\Omega}
\newcommand{\quext}{\quad\text}

\DeclareMathOperator{\deriv}{d}

\newcommand{\dit}{\deriv\!t}
\newcommand{\dis}{\deriv\!s}
\newcommand{\dix}{\deriv\!x}

\newcommand{\perogni}{\forall}

\newcommand{\dual}[1]{(\!({#1})\!)}
\newcommand{\duav}[1]{\langle\!\langle{#1}\rangle\!\rangle}

\numberwithin{equation}{section}
\mathchardef\emptyset="001F

\newtheorem{theorem}{Theorem}[section]

\newtheorem{lemma}[theorem]{Lemma}

\newtheorem{definition}[theorem]{Definition}
\theoremstyle{definition}

\numberwithin{equation}{section}
\begin{document}

\title{On the viscous Cahn-Hilliard equation with singular 
 potential and inertial term}

\author{%
Riccardo Scala\\
Dipartimento di Matematica, Universit\`a di Pavia,\\
Via Ferrata~1, 27100 Pavia, Italy\\
E-mail: {\tt riccardo.scala@unipv.it}
\and 
Giulio Schimperna\\
Dipartimento di Matematica, Universit\`a di Pavia,\\
Via Ferrata~1, 27100 Pavia, Italy\\
E-mail: {\tt giusch04@unipv.it}
}


\maketitle
\begin{abstract}
 We consider a relaxation of the viscous Cahn-Hilliard equation
 induced by the second-order inertial term~$u_{tt}$.  
 The equation also contains a semilinear term $f(u)$ 
 of ``singular'' type. Namely, the function $f$ is defined only on 
 a bounded interval of $\RR$ corresponding to the physically admissible 
 values of the unknown $u$, and diverges as $u$ approaches 
 the extrema of that interval. In view of its interaction 
 with the inertial term $u_{tt}$, the term $f(u)$
 is difficult to be treated mathematically.
 Based on an approach originally devised for the strongly damped wave equation,
 we propose a suitable concept of weak solution based on duality
 methods and prove an existence result.
\end{abstract}

\noindent {\bf Key words:}~~Cahn-Hilliard equation, inertia, 
 weak formulation, maximal monotone operator, duality.

\vspace{2mm}

\noindent {\bf AMS (MOS) subject clas\-si\-fi\-ca\-tion:}~~35K67, 35L85,
 46A20, 47H05, 80A22.


\section{Introduction}
\label{sec:intro}

The celebrated Cahn-Hilliard equation was proposed to describe phase separation
phenomena in binary systems \cite{CH}. Its ``standard'' version has the form
of a semilinear parabolic fourth order equation, namely
\begin{equation}\label{CH-par}
  u_t - \Delta( -\Delta u + f(u) )=0.
\end{equation}
Here the unknown $u$ stands for the relative 
concentration of one phase, or component,
in a binary material, and $f$ is the derivative of a non-convex potential $F$
whose minima represent the energetically more favorable configurations
usually attained in correspondance, or in proximity, 
of pure phases or concentrations.
In view of the fact that $u$ is an order parameter, often
it is normalized in such a way that the pure states
correspond to the values $u=\pm 1$, whereas $-1<u<1$ denotes the
(local) presence of a mixture. We will also adopt this convention. In such
setting the values $u\not\in [-1,1]$ are generally interpreted as
``nonphysical'' and should be somehow excluded. In view
of the fourth-order character of~\eqref{CH-par}, no 
maximum principle is available for $u$. Hence, the {\sl constraint}\/
$u\in [-1,1]$ is generally enforced by assuming $F$ to 
be defined only for $u\in(-1,1)$ (or for $u\in[-1,1]$; both 
choices are admissible under proper structure 
conditions) and to be identically $+\infty$ outside
the interval $[-1,1]$. A relevant example is given 
by the so-called {\sl logarithmic potential}
\begin{equation}\label{logpot}
  F(u)=(1-u)\log(1-u) + (1+u)\log(1+u) - \frac{\lambda}2 u^2,
   \quad \lambda \ge 0,
\end{equation}
where the last term may induce nonconvexity.
Such a kind of potential is generally termed as 
a {\sl singular}\/ one and its occurrence may give rise to mathematical
difficulties in the analysis of the system. For this reason, 
singular potentials are often replaced by ``smooth'' approximations
like the so-called {\sl double-well}\/ potential taking, 
after normalization, the form $F(u)= (1-u^2)^2$. Of course,
in the presence of a smooth double-well potential, solutions
are no longer expected to satisfy the physical constraint $u\in[-1,1]$.

The mathematical literature devoted to~\eqref{CH-par}
is huge and the main properties of the solutions 
in terms of regularity, qualitative behavior, and asymptotics 
are now well-understood, also in presence of singular potentials 
like~\eqref{logpot} (cf., e.g., \cite{N-C1,MZ1} and the 
references therein). Actually, in recent years, the attention has moved to
more sophisticated versions of~\eqref{CH-par} related to specific
physical situations. Among these, we are interested here 
in the so-called {\sl hyperbolic relaxation}\/ of the equation. 
This can be written as
\begin{equation}\label{CH-hyp}
  \alpha u_{tt} + u_t - \Delta( -\Delta u + f(u) ) = 0,
\end{equation}
where $\alpha>0$ is a (small) relaxation parameter and the new
term accounts for the occurrence of ``inertial'' effects.
Equation~\eqref{CH-hyp} may be used in order to describe 
strongly non-equilibrium decomposition generated by
deep supercooling into the spinodal region occurring in certain materials 
(e.g., glasses), see \cite{Gal,GJ}. From the mathematical point of view, 
equation~\eqref{CH-hyp} carries many similarities with the semilinear
(damped) wave equation, but is, however, much more delicate to deal with.
For instance, in space dimension $N=3$ the existence of global in time
strong solutions is, up to our knowledge, an open issue also in
the case when $f$ is a globally Lipschitz (nonlinear)
function \cite{GSSZ}, 
whereas for $N=2$ the occurrence of a critical exponent
is observed in case $f$ has a polynomial growth \cite{GSZ1,GSZ2}. 
The situation is somehow more satisfactory
in space dimension $N=1$ (cf., e.g., \cite{ZM1,ZM2}) due to 
better Sobolev embeddings (in particular all solutions taking values
in the ``energy space'' are also uniformly bounded).
It is however worth noting that, in the case when $f$ is singular,
even the existence of (global) {\sl weak}\/ solutions 
is a mathematically very challenging problem. Indeed, at least 
up to our knowledge, this seems to be an open issue even in
one space dimension.

The picture is only partially more satisfactory when one considers 
a further relaxation of the equation containing a ``strong damping''
(or ``viscosity'') term, namely
\begin{equation}\label{CH-strong}
  \alpha u_{tt} + u_t - \Delta( \delta u_t -\Delta u + f(u) ) = 0,
\end{equation}
with $\delta>0$ (a physical justification for this equation
is given, e.g., in \cite{N-C}). The new term induces additional regularity 
and some parabolic smoothing effects, and, for this reason,
\eqref{CH-strong} is mathematically more tractable
in comparison to~\eqref{CH-hyp}. Indeed, existence, regularity and large time
behavior of solutions have been analyzed in a number of papers
(cf., e.g., \cite{Bonfoh,BGM,GGMP1,GPS,Kan} and references therein).
In all these contributions, however, $f$ is taken as a smooth 
function of at most polynomial growth at infinity.
Here, instead, we will consider \eqref{CH-strong} with the choice of
a singular function $f$.

To explain the related difficulties, the main 
point stands, of course, in the low
number of available a-priori estimates. This is 
a general feature of equations of the second order in time, and,
as a consequence, approximating sequences satisfy very poor compactness
properties. In particular, the second order term $u_{tt}$
can be only controlled in a space like $L^1(0,T;X)$, 
where $X$ is a Sobolev space of negative order. In view of the 
bad topological properties of $L^1$, this implies
that in the limit the term $u_t$ cannot be shown to be (and, in fact,
is not expected to be) continuous in time, but
only of bounded variation. In particular, it may present 
{\sl jumps}\/ with respect to the time variable. In turn, 
the occurrence of these jumps is strictly connected to the fact that it
is no longer possible to compute the singular term $f(u)$ 
in the ``pointwise'' sense.  

Indeed, in the weak formulation $f(u)$ is 
suitably reinterpreted in the distributional sense, and, in
particular, concentration phenomena may occur.
This idea comes from the theory of convex integrals
in Sobolev spaces introduced in the celebrated paper 
by Brezis~\cite{brezisart}
and later developed and adapted to cover a number of 
different situations (cf., e.g., \cite{BCGG,BGRT,SP} 
and references therein). In our former paper
in collaboration with E.~Bonetti and E.~Rocca \cite{BRSS} 
we have shown that this method can be adapted to treat
equations of the second order in time. Actually, 
using duality methods in Sobolev spaces of {\sl parabolic type}\/
(i.e., depending both on space and on time variables), 
we may provide the required relaxation of the term $f(u)$
accounting for the possible occurrence of 
concentration phenomena with respect to time.
The reader is referred to~\cite{BRSS} for further considerations 
and extended comments and examples.

Equation \eqref{CH-strong} will be considered here 
in the simplest mathematical setting. Namely, we will settle it 
in a smooth bounded domain $\Omega\subset \RR^N$, 
$N\le 3$ (we remark however that the results could be easily extended
to any spatial dimension), in a fixed reference interval 
$(0,T)$ of arbitrary length, and with homogeneous
Dirichlet boundary conditions. Then, existence of weak solutions
will be proved by suitably adapting the approach of~\cite{BRSS}.
It is worth observing that, as happens for the mentioned strongly 
damped wave equation and for other similar models, an alternative 
weak formulation could be given by restating the problem in the
form of a variational inequality. However, as noted in \cite{BRSS},
we believe the concept of solution provided here to be somehow
more flexible. In particular, with this method we may provide an
explicit characterization of the (relaxed) term $f(u)$
(which may be thought as a physical quantity representing 
the vincular reaction provided by the constraint) in
terms of regularity (for instance, for equation~\eqref{CH-strong}
concentration phenomena are expected to occur
{\sl only}\/ with respect to the time variable $t$).
Moreover, we can prove that at least {\sl some}\/ weak 
solutions satisfy a suitable form of the energy 
inequality. This can be seen as a sort of selection principle
for ``physical'' solutions (note, indeed, that uniqueness is
not expected to hold).

\smallskip

The plan of the paper is the following: in the next Section~\ref{sec:main} we 
introduce our assumptions on coefficients and data and state our main result
regarding existence of at least one solution to a suitable weak formulation
of equation~\eqref{CH-strong}. The proof of this theorem is then carried out in 
Section~\ref{sec:proof} by means of an approximation -- a priori estimates --
compactness argument.


\section{Main result}
\label{sec:main}


\subsection{Preliminaries}
\label{subsec:prel}

We consider the viscous Cahn-Hilliard equation with inertia:
\begin{equation}\label{CH}
 \begin{cases}
  \alpha u_{tt}+u_t-\Delta w=0,\\
  w=\delta u_t-\Delta u+\beta(u)-\lambda u.
 \end{cases}
\end{equation}
Here the coefficients $\alpha$ and $\delta$ are strictly positive constants,
whereas $\lambda\ge 0$. Moreover, $\beta$ is a maximal monotone operator 
in $\R\times\R$ satisfying 
\begin{equation}\label{hyp:beta}
  \overline{D(\beta)}=[-1,1], \quad 0\in\beta(0).
\end{equation}
Actually, $\beta$ represents the monotone part of $f(u)$ 
(cf.~\eqref{CH-hyp}). The {\sl domain}~$D(\beta)$ has been normalized
just for mathematical convenience. Following \cite{Ba}, 
there exists a convex and lower semicontinuous 
function $j:\R\rightarrow[0,+\infty]$ such that $\beta=\partial j$, 
$\overline{D(j)}=[-1,1]$, and $j(0)=\min j=0$.
For all $\eps\in(0,1)$ we denote by $j^\eps:\R\rightarrow[0,+\infty)$ the 
{\sl Moreau-Yosida regularization}\/ of $j$, 
and by $\beta^\eps:=\partial j^\eps=(j^\eps)'$ 
the corresponding {\sl Yosida approximation}\/ of $\beta=\partial j$. 

By a direct check (cf.~also \cite[Appendix~A]{MZ1}), one may prove, based on 
\eqref{hyp:beta}, that there exist constants 
$c_1>0$ and $c_2\geq0$ independent of $\eps$, such that 
\begin{equation}\label{bound:l1}
  \beta^\eps(r)r \ge c_1 |\beta^\eps(r)| - c_2.
\end{equation}
Let us also introduce some functional spaces: we set $H:=L^2(\Om)$ and $V:=H^1_0(\Om)$, 
so that $V'=H^{-1}(\Om)$. Moreover, we put
$$
  \mathcal V:=H^1(0,T;H),
$$
and, for all $t\in(0,T]$,
$$
  \mathcal V_t:=H^1(0,t;H).
$$
We denote by $(\cdot,\cdot)$ and  $\langle\cdot,\cdot\rangle$ the scalar 
product in $H$ and the duality pairing between $V'$ and $V$, respectively. 
The scalar products on $L^2(0,T;H)$ and on $L^2(0,t;H)$, 
for $t\in(0,T)$, are indicated respectively by
$$
  \dual{\cdot,\cdot}\;\;\;\;\text{ and by }\;\;\;\;\dual{\cdot,\cdot}_t.
$$
Correspondingly, the duality products between $\mathcal{V}$ 
and $\mathcal V'$ and between $\mathcal{V}_t$ and $\mathcal V'_t$ are noted as
$$
  \duav{\cdot,\cdot}\;\;\;\;\text{ and }\;\;\;\;\duav{\cdot,\cdot}_t,
$$
respectively.

Next, we indicate by $A:D(A)\rightarrow H$, with domain $D(A):=H^2(\Om)\cap H^1_0(\Om)$,
the Laplace operator with homogeneous Dirichlet boundary condition seen as 
an unbounded linear operator on~$H$. Hence, $A$ is strictly positive and
its powers $A^s$ are well defined for all $s\in\R$. In particular,
$D({A^{1/2}})=H^1_0(\Om)=V$. Moreover, $A$ may be extended to the 
space $V$ and it turns out that $A:V\rightarrow V'$ is an isomorphism.
In particular, $V'$ is a Hilbert space when endowed with the scalar product
$$ 
  (u,v)_*:= \langle v,A^{-1}u \rangle = \langle u,A^{-1}v \rangle 
   \quext{for }\,u,v\in V'.
$$
The associated norm is then given by $\|u\|_{V'}^2=(u,u)_*$ 
for $u\in V'$. Correspondingly,
the scalar products of the spaces $L^2(0,T;V')$ and $L^2(0,t;V')$ are denoted by 
$$ 
  (\!(\cdot,\cdot)\!)_*\;\;\;\;\;\;(\!(\cdot,\cdot)\!)_{*,t}
$$
respectively. In particular, we have
$$ 
  (\!(u,v)\!)_* = \int_0^T \langle v,A^{-1}u \rangle \dit 
   \quext{for }\,u,v\in L^2(0,T;V'),
$$
with a similar characterization holding 
for~$(\!(\cdot,\cdot)\!)_{*,t}$.


\subsection{Relaxation of the constraint}
\label{subsec:constr}

We now provide a brief sketch of the relaxation of $\beta$ mentioned
in the introduction, referring to~\cite[Sec.~2]{BRSS} for additional
details. 
%
%
First of all, we introduce the functional $J:H\to [0,+\infty]$,
$J(u):=\int_{\Om}j(u)\,\dix$ for all $u\in H$, whose 
value is intended to be $+\infty$ if $j(u)\notin L^1(\Om)$. 
Moreover it is convenient to define
\begin{equation}
  \mathcal J(u):=\int_0^T\int_{\Om}j(u)\,\dix\,\dit
   \;\;\;\;\;\forall\, u\in L^2((0,T)\times\Om),
\end{equation}
and its counterpart on $(0,t)$, namely
\begin{equation}
  \mathcal J_t(u):=\int_0^t\int_{\Om}j(u)\,\dix\,\dis
   \;\;\;\;\;\forall\, u\in L^2((0,t)\times\Om).
\end{equation}
Then, the relaxed version of $\beta$ will be intended as a
maximal monotone operator in the duality couple $\calV \times \calV'$.
Indeed, we first introduce $\mathcal J_{\mathcal V}:=\mathcal J\llcorner_{\mathcal V}$, 
the restriction of $\mathcal J$ to $\mathcal V$. Then, we consider its 
subdifferential $\partial\mathcal J_{\mathcal V}$ with respect 
to the duality pairing between $\mathcal V$ and $\mathcal V'$. 
Namely, for $\xi\in \mathcal V'$ and $u\in\mathcal V$, we say that 
\begin{align}
   \xi\in \partial \mathcal J_{\mathcal V}(u)\;\;\Longleftrightarrow\;\;
  \mathcal J_{\mathcal V}(z) \geq\duavw{\xi,z-u}+\mathcal J_{\mathcal V}(u)\;\;\;\;
  \forall \,z\in\mathcal V.
\end{align}
In order to emphasize that $\partial \mathcal J_{\mathcal V}$ 
consists in a relaxation of~$\beta$, we will simply note
$\partial \mathcal J_{\mathcal V}=:\beta_w$ ($w$ standing for ``weak'').
Proceeding in a similar way for the functional $\mathcal J_t$, we define the 
subdifferential $\partial\mathcal J_{t,\mathcal V_t}$ of the operator 
$\mathcal J_{t,\mathcal V_t}:=\mathcal J_t\llcorner_{\mathcal V_t}$. This 
will be indicated simply by $\beta_{w,t}$.

In this setting it is not true anymore that an {\sl element}\/ $\xi$ of the
{\sl set}\/ $\beta_w(u)$ (recall that $\beta$ is a multivalued operator and, 
as a consequence, $\beta_w$ may
be multivalued as well) admits a ``pointwise'' interpretation
as ``$\xi(t,x)=\beta(u(t,x))$''.
Indeed, $\xi$ belongs to the negative order Sobolev space $\calV'$
and concentration phenomena are expected to occur. Nevertheless, the 
maps $\beta^\eps$ still provide a suitable approximation of $\beta_w$.
Referring the reader to~\cite{BRSS,S15} for additional details and comments,
we just mention here some basic facts. First of all, let us 
define $J^\eps(u):=\int_{\Om}j^\eps(u)\,\dix$ and 
$\calJ^\eps(u):=\int_0^T\int_{\Om}j^\eps(u)\,\dix\,\dit$.
Then, one may prove that the functionals $\mathcal J^\eps$ converge 
to~$\mathcal J$ in the sense of Mosco-convergence 
with respect to the topology of $L^2(0,T;H)$. Moreover, their restrictions
to $\calV$ Mosco-converge to $\mathcal J_{\calV}$ 
in the topology of $\calV$. The analogue of these properties 
also holds for restrictions to time subintervals $(0,t)$.
Referring the reader to~\cite[Chap.~3]{Att}
for the definition and basic properties of Mosco-convergence, here
we just recall that this convergence notion for functionals
implies (and is in fact equivalent to) a related notion
of convergence for their subdifferentials, called graph-convergence
(or G-convergence). Namely, noting that the function $\beta^\eps$ represents 
the subdifferential of $\mathcal J^\eps$ both with respect to 
the topology of $L^2(0,T;H)$ and to that of $\mathcal V$,
it turns out that the operators $\beta^\eps$, if identified
with their graphs, G-converge to $\beta$ in the topology
of $L^2(0,T;H)\times L^2(0,T;H)$ and G-converge to $\beta_w$ in the topology
of $\calV\times\calV'$. As a consequence of the latter property, 
we may apply the so-called 
%
%
{\sl Minty's trick}\/ in the duality between $\calV$ and $\calV'$.
This argument will be the main tool we will use in order to take the limit
in the approximation of the problem and can be simply stated in this way:
once one deals with a sequence $\{v^\eps\}\subset \calV$
satisfying $v_\eps\rightharpoonup v$ weakly in $\mathcal V$ and 
$\beta^\eps(v_\eps)\rightharpoonup \xi$ weakly in $\mathcal V'$, 
then the inequality  
\begin{equation}\label{minty}
   \limsup_{\eps\searrow0}\, \duav{\xi_\eps,v_\eps}\leq \duav{\xi,v}
\end{equation}
implies that $\xi\in \beta_w(v)$. In other words, $\xi$ is identified as an
element of the set $\beta_w(v)\subset \calV'$.


\subsection{Statement of the main result}
\label{subsec:res}

We start with presenting our basic concept of weak solution,
which can be seen as an adaptation of~\cite[Def.~2.2]{BRSS}.
\begin{definition}\label{def:sol}
 A couple $(u,\eta)$ is called a weak solution to the initial-boundary value problem for
 the viscous Cahn-Hilliard equation with inertia whenever the following 
 conditions hold:\\[2mm]
 {\sl (a)}~There hold the regularity properties
 \begin{align}\label{reg:ut}
   & u_t \in BV(0,T;H^{-4}(\Om)) \cap L^{\infty}(0, T ;V') \cap L^2(0, T ;H),\\
  \label{reg:u}
   & u \in L^{\infty}(0, T ;V) \cap L^2(0, T ;D(A)),\\ 
  \label{reg:xi}
   & \eta \in \calV'.
 \end{align}
 {\sl (b)}~For any test function $\varphi\in \mathcal V$,
 there holds the following weak version of~\eqref{CH}:
 \begin{align}\nonumber
   & \alpha(u_{t}(T),\varphi(T))_*
    - \alpha (u_1,\varphi(0))_*
    - \alpha\dual{u_{t},\varphi_t}_*
    + \dual{u_t,\varphi}_*\\
  \label{eq:w}
   & \mbox{}~~~~~ 
    + \delta\dual{u_t,\varphi}
    + \dual{A^{1/2} u, A^{1/2}\varphi}
    + \duav{\eta,\varphi}
    - \lambda\dual{u,\varphi}=0.  
 \end{align}
 Moreover, for all $t\in [0,T]$ there exists 
 $\eta_{(t)}\in \mathcal{V}'$ such that 
 \begin{align}\nonumber
  & \alpha (u_{t}(t),\varphi(t))_*
   - \alpha (u_1,\varphi(0))_*
   - \alpha\dual{u_{t},\varphi_t}_{*,t}
   + \dual{u_t,\varphi}_{*,t}\\
 \label{eq:wt}
   & \mbox{}~~~~~
   + \delta\dual{u_t,\varphi}_t
   + \dual{A^{1/2} u, A^{1/2}\varphi}_t
   + \duav{\eta_{(t)},\varphi}_t
   - \lambda\dual{u,\varphi}_t=0,
 \end{align}
 for all $\varphi\in \mathcal V_t$.\\[2mm]
 {\sl (c)}~The functionals $\eta$ and $\eta_{(t)}$ satisfy
 \begin{equation}\label{graph}
   \eta\in\beta_w(u), \qquad \eta_{(t)}\in \beta_{w,t}(u\llcorner_{(0,t)})
    \;\;\text{ for all }t\in(0,T),
 \end{equation}
 and the following compatibility condition holds true:
 \begin{equation}\label{compa}
   \duav{\eta_{(t)},\varphi}_t = \duav{\eta,\bar\varphi}
    \;\;\text{ for all }
    \varphi\in\mathcal V_{t,0}\;\;\text{ and all }t\in[0,T),
 \end{equation}
 where $\calV_{t,0}:=\{ \varphi\in \calV_t:~\varphi(t)=0\}$
 and $\bar\varphi$ is the trivial extension of 
 $\varphi\in \mathcal V_{t,0}$ to $\mathcal V$, i.e.,
 $\bar\varphi(s)=\varphi(t)=0$ for all $s\in (t,T]$.\\[2mm]
 {\sl (d)}~There holds the Cauchy condition
 \begin{equation}\label{init}
   u|_{t=0} = u_0 \quext{a.e.~in }\,\Omega.
 \end{equation}
\end{definition}
\noindent%
Correspondingly, we conclude this section with our main result, 
stating existence of at least one weak solution.
\begin{theorem}\label{main}
 Let $T>0$ and let the initial data satisfy 
 \begin{equation}\label{hp:init}
   u_0\in V, \quad j(u_0)\in L^1(\Omega), 
    \qquad u_1\in H.
 \end{equation}
 Then, there exists a solution $(u,\eta)$
 to the viscous Cahn-Hilliard equation with inertia in the 
 sense of\/ {\rm Def.~\ref{def:sol}}. Moreover, 
 $u$ satisfies the\/ {\rm energy inequality}
 \begin{align}\label{energ:eq}
   & \frac{\alpha}{2} \|u_t(t_2)\|_{V'}^2
    + \frac{1}{2}\|A^{1/2} u(t_2)\|^2_{H} 
    + J(u(t_2))
    - \frac{\lambda}{2}\|u(t_2)\|^2_H
    + \int_{t_1}^{t_2}\big( \delta\|u_t\|_H^2
    + \|u_t\|^2_{V'}\big) \,\dis\nonumber\\
   & \mbox{}~~~~~
    \le \frac{\alpha}{2}\|u_t(t_1)\|_{V'}^2
    + \frac{1}{2}\|A^{1/2} u(t_1)\|^2_{H}
    + J(u(t_1))
    - \frac{\lambda}{2}\|u(t_1)\|^2_H,
 \end{align}
 for\/ {\rm almost every} $t_1\in [0, T)$ (surely including $t_1 = 0$)
 and\/ {\rm every} $t_2 \in (t_1,T]$. 
\end{theorem}


\section{Proof of Theorem~\ref{main}}
\label{sec:proof}


\subsection{Approximation}
\label{subsec:appro}

We consider a regularization of system~\eqref{CH}, 
namely for $\eps\in(0,1)$ we denote by $(u^\eps,w^\eps)$ the solution to 
\begin{align}\label{CH:eps:1}
  & \alpha u^\eps_{tt}+u^\eps_t+A w^\eps = 0,\\
 \label{CH:eps:2}
  & w^\eps=\delta u^\eps_t+A u^\eps+\beta^\eps(u^\eps)-\lambda u^\eps,
\end{align}
coupled with the initial conditions
\begin{equation}\label{init:eps}
   u^\eps|_{t=0} = u_0^\eps ~~\text{and }\,\,
   u^\eps_t|_{t=0} = u_1^\eps, 
     \quext{a.e.~in }\,\Omega.
\end{equation}
Recall that $\beta^\eps$ was defined in Subsec.~\ref{subsec:prel}.
The following result provides existence of a unique smooth solution 
to \eqref{CH:eps:1}-\eqref{init:eps} 
once the initial data are suitably regularized:
\begin{theorem}\label{exist:eps}
 Let $T>0$, $u^\eps_0\in D(A)=H^2(\Om)\cap V$, 
 $u^\eps_1\in D(A^{1/2})=V$. Then there exists a unique
 function $u^\eps$ with 
 \begin{align}\label{regee1}
  &u^\eps\in  W^{1,\infty}(0,T;H)\cap H^1(0,T;V)\cap L^\infty (0,T;D(A)),\\
  \label{regee2}
  &u^\eps_t\in W^{1,\infty}(0,T;D({A^{-1}})), 
 \end{align}
 satisfying \eqref{CH:eps:1}-\eqref{init:eps}.
 Moreover, for every $t_1,t_2\in[0,T]$, there holds the approximate energy balance
 \begin{align}\label{energ:eps}
   & \frac{\alpha}{2}\|u_t^\eps(t_2)\|_{V'}^2
    + \frac{1}{2}\|A^{1/2} u^\eps(t_2)\|^2_{H}
    + J^\eps(u^\eps(t_2)) 
    - \frac{\lambda}{2}\|u^\eps(t_2)\|^2_H
    + \int_{t_1}^{t_2} 
      \big( \delta\|u_t^\eps\|_H^2 + \|u_t^\eps\|^2_{V'} \big)\,\dis \nonumber\\
  & \mbox{}~~~~~
    = \frac{\alpha}{2}\|u_t^\eps(t_1)\|_{V'}^2
    + \frac{1}{2}\|A^{1/2} u^\eps(t_1)\|^2_{H}
    + J^\eps(u^\eps(t_1)) 
    - \frac{\lambda}{2}\|u^\eps(t_1)\|^2_H.
 \end{align}
\end{theorem}
\noindent%
The proof of the above result is standard (see, e.g., \cite[Thm.~2.1]{GGMP2}).
Actually, one can replicate the a-priori estimates corresponding to
the regularity properties \eqref{regee1}-\eqref{regee2}
by multiplying \eqref{CH:eps:1} by $u^\eps_t$, \eqref{CH:eps:2} by $A u^\eps_t$,
and using the Lipschitz continuity of $\beta^\eps$.
The regularity of $\beta^\eps$ is also essential for having
uniqueness, as one can show via standard contractive methods.
Then, to prove the energy equality it is sufficient to test \eqref{CH:eps:1} 
by $A^{-1} u^\eps_t$, \eqref{CH:eps:2} by $u^\eps_{t}$, and 
integrate the results with respect to the time and space variables.
It is worth observing that these test functions are {\sl admissible}\/
thanks to the regularity properties \eqref{regee1}-\eqref{regee2}. As a consequence
of this fact, we can apply standard chain-rule formulas
to obtain that \eqref{energ:eps} holds with the {\sl equal}\/ sign,
which will no longer be the case in the limit.

\smallskip

As a first step in the proof of Theorem~\ref{main}, we need to 
specify the required regularization of the initial data:
\begin{lemma}\label{iniz_data}
 Let~\eqref{hp:init} hold. Then there exist two families
 $\{u^\eps_0\}\subset D(A)\cap V$
 and $\{u_1^\eps\}\subset V$, $\eps\in(0,1)$, satisfying
 \begin{align}\label{apprin:1} 
   & J^\eps(u^\eps_0)\leq J(u_0)~~\perogni\eps>0\;\;\;\;\;\;
    \text{ and }\;\;u^\eps_0\rightarrow u_0~~\text{in }\,V,\\
   \label{apprin:2} 
   & u^\eps_1\rightarrow u_1~~\text{in }\,H.
 \end{align}
\end{lemma}
\noindent%
Also the above lemma is standard. Indeed, one can construct
$u^\eps_0$, $u^\eps_1$ by simple singular perturbation methods 
(see, e.g., \cite[Sec.~3]{SP}). 
Let us then consider the solutions $u^\eps$ to the regularized system 
\eqref{CH:eps:1}-\eqref{init:eps}
with the initial data provided by Lemma~\ref{iniz_data}.
Then, taking a test function $\varphi\in\mathcal V$, multiplying 
\eqref{CH:eps:1} by $A^{-1}\varphi$, \eqref{CH:eps:2} by $\varphi$, 
and performing standard manipulations, one 
can see that $u^\eps$ also satisfies the weak formulation 
(compare with~\eqref{eq:w})
\begin{align}\nonumber
  & \alpha(u^\eps_{t}(T),\varphi(T))_*
   - \alpha (u^\eps_1,\varphi(0))_*
   - \alpha\dual{u^\eps_{t},\varphi_t}_*
   + \dual{u^\eps_t,\varphi}_*\\
 \label{eq:wee}
  & \mbox{}~~~~~
   + \delta\dual{u^\eps_t,\varphi}
   + \dual{A^{1/2} u^\eps, A^{1/2}\varphi}
   + \dual{\beta^\eps(u^\eps),\varphi}
   - \lambda\dual{u^\eps,\varphi}=0.  
\end{align}
Correspondingly, the analogue over subintervals $(0,t)$ also holds. Namely,
for $\varphi\in\mathcal V_t$ one has (compare with~\eqref{eq:wt})
\begin{align}\nonumber
  & \alpha(u^\eps_{t}(t),\varphi(t))_*
   - \alpha (u^\eps_1,\varphi(0))_*
   - \alpha\dual{u^\eps_{t},\varphi_t}_{*,t}
   + \dual{u^\eps_t,\varphi}_{*,t}\\
 \label{eq:weet}
  & \mbox{}~~~~~
   + \delta\dual{u^\eps_t,\varphi}_t
   + \dual{A^{1/2} u^\eps, A^{1/2}\varphi}_t
   + \dual{\beta^\eps(u^\eps),\varphi}_t
   - \lambda\dual{u^\eps,\varphi}_t=0.     
\end{align}


\subsection{A priori estimates}
\label{subsec:apriori}

We now establish some a-priori estimates for $u^\eps$. The estimates
will be uniform in $\eps$ and permit us to take $\eps\searrow 0$ 
at the end. First of all, the energy balance \eqref{energ:eps} 
and the uniform bounded properties \eqref{apprin:1}-\eqref{apprin:2} 
of approximating initial data
provide the existence of a constant $M>0$, independent of $\eps$, 
such that the following bounds hold true:
\begin{subequations}\label{ests}
  \begin{align}
  &\|u^\eps\|_{L^{\infty}(0,T;V)}\leq M\label{est1},\\
  &\|u^\eps\|_{H^1(0,T;V')}\leq M,\label{est1b}\\
  &\delta^{1/2}\| u^\eps\|_{H^1(0,T;H)}\leq M,\label{est2}\\
  &\alpha^{1/2}\| u^\eps\|_{W^{1,\infty}(0,T;V')}\leq M,\label{est2b}\\
  &\| j^\eps(u^\eps)\|_{L^\infty(0,T;L^1(\Om))}\leq M,
\end{align}
\end{subequations}
for all $\eps\in(0,1)$. More precisely, thanks to the fact that,
for every (fixed) $\eps\in (0,1)$, $u^\eps_t$ lies in $C^0([0,T];V')$
by~\eqref{regee1}-\eqref{regee2},
we are allowed to evaluate $u^\eps_t$ pointwise in time.
Hence, \eqref{est2b} may be complemented by
\begin{equation}\label{est2c}
  \| u^\eps_t(t) \|_{V'}\leq M \;\;\;\;\;
   \text{ for every }\,t\in[0,T],
\end{equation}
and in particular, for $t=T$. Analogously, 
thanks to $u^\eps \in C^0([0,T];V)$,
in addition to~\eqref{est1} we also have
\begin{align}
  & \|u^\eps(t)\|_{V}\leq M\;\;\;\text{ for every }\,t\in [0,T].\label{est1*}
\end{align}
Next, taking $\varphi=u^\eps$ in \eqref{eq:wee} and rearranging terms, we infer
\begin{align}\no
 & \int_0^T\int_\Om\beta^\eps(u^\eps) u^\eps\,\dix\,\dit
   \leq \alpha\|u^\eps_t(T)\|_{V'}\|u^\eps(T)\|_{V'}
   +\alpha\|u^\eps_1\|_{V'}\|u^\eps_0\|_{V'}\\
 \no
  & \mbox{}~~~~~
   + \alpha\|u^\eps_t\|_{L^2(0,T;V')}^2 
   + \|u^\eps_t\|_{L^2(0,T;V')}\|u^\eps\|_{L^2(0,T;V')} \\
 \label{estxx}
  & \mbox{}~~~~~
   + \delta \|u^\eps_t\|_{L^2(0,T;H)}\|u^\eps\|_{L^2(0,T;H)}
   + \|A^{1/2} u^\eps\|^2_{L^2(0,T;H)}
   + \lambda\|u^\eps\|_{L^2(0,T;H)}^2.
\end{align}
Then, thanks to estimates \eqref{ests},  \eqref{est2c} and \eqref{est1*},
we may check that the right hand side of~\eqref{estxx}
is bounded uniformly with respect
to~$\eps$. Consequently, using also \eqref{bound:l1}, we infer
\begin{equation}\label{estbeta}
  \|\beta^\eps(u^\eps)\|_{L^1(0,T;L^1(\Om))}\leq M.
\end{equation}
Now, since we assumed $N\leq3$, we know that $L^1(\Om)\subset D(A^{-1})$,
the latter being a closed subspace of $H^{-2}(\Om)$.
Moreover, $A$ can be extended to a bounded linear operator
$A:D(A^{-1})\rightarrow D(A^{-2}) \subset H^{-4}(\Om)$. Then, letting  
$X:=H^{-4}(\Om)$ (note that for $N>3$ the argument still works up to
suitably modifying the choice of~$X$) and rewriting 
\eqref{CH:eps:1}-\eqref{CH:eps:2} as a single equation, i.e.,
\begin{equation}\label{CH:eps:}
  \alpha u^\eps_{tt}+u^\eps_t
   + \delta A u^\eps_t+A^2 u^\eps
   + A (\beta^\eps(u^\eps)) -\lambda A u^\eps = 0,
\end{equation}
we may check by a comparison of terms that
\begin{equation}\label{est5}
 \alpha \|u^\eps_t\|_{W^{1,1}(0,T;X)}\leq M.
\end{equation}
Actually, we used here the 
estimates \eqref{ests} together with \eqref{estbeta}.

Next, thanks to the last of~\eqref{regee1}, we are allowed
to multiply \eqref{CH:eps:1} by $u^\eps$ 
and \eqref{CH:eps:2} by $Au^\eps$. Using the monotonicity
of $\beta^\eps$ and the bounds \eqref{ests}, standard 
arguments lead us to the additional estimate
\begin{equation}\label{est6}
  \| u^\eps\|_{L^2(0,T;D(A))}\leq M,
\end{equation}
still holding for $M>0$ independent of $\eps$.

Finally, for all $\varphi\in\mathcal V_t$
we can compute from \eqref{eq:weet}
\begin{align}
  & \left|\int_0^t\langle\beta^\eps(u^\eps),\varphi\rangle\,\dis \right|
   \leq \alpha\|u^\eps_t(t)\|_{V'}\|\varphi(t)\|_{V'}
     + \alpha \|u^\eps_1\|_{V'}\|\varphi(0)\|_{V'}
     + \alpha\|u^\eps_t\|_{L^2(0,t;V')}\|\varphi_t\|_{L^2(0,t;V')}\nonumber\\
  & \mbox{} ~~~~~
    + \|u^\eps_t\|_{L^2(0,t;V')}\|\varphi\|_{L^2(0,t;V')}
    + \delta\|u^\eps_t\|_{L^2(0,t;H)}\|\varphi\|_{L^2(0,t;H)}\nonumber\\
 \label{est6b}
  & \mbox{} ~~~~~
   + \| u^\eps\|_{L^2(0,t;D(A))}\|\varphi\|_{L^2(0,t;H)}
   + \lambda\|u^\eps\|_{L^2(0,t;H)}\|\varphi\|_{L^2(0,t;H)}
\end{align}
and the right-hand side, by \eqref{ests}, \eqref{est2c} and \eqref{est6}, 
is less or equal than $C\|\varphi\|_{\mathcal V_t}$,
with $C$ depending only on the (controlled) norms of $u^\eps$.
Hence it follows that there exists a constant $M>0$ independent
of $\eps$ such that
\begin{align}
  \|\beta^\eps(u^\eps)\|_{\mathcal V_t'}\leq M,\label{est:beta:V'}
\end{align}
for every $t\in(0,T]$. In particular, $\|\beta^\eps(u^\eps)\|_{\mathcal V'}\leq M$.


\subsection{Passage to the limit}
\label{subsec:limit}

Using the estimates obtained above, we now aim to pass to the 
limit as $\eps\searrow 0$ in the weak formulation \eqref{eq:wee}.
Firstly, \eqref{ests}, \eqref{est6} and \eqref{est:beta:V'} imply that there exist
$u\in W^{1,\infty}(0,T;V')\cap H^1(0,T;H)\cap L^\infty(0,T;V)\cap L^2(0,T;D(A))$ and 
$\eta\in \mathcal V'$ such that
\begin{subequations}\label{conv}
 \begin{align}
  &u^\eps\rightharpoonup u \;\;\;\text{ weakly star in } W^{1,\infty}(0,T;V')
  \text{ and weakly in }L^2(0,T;D(A)),\label{conv1}\\
    &u^\eps\rightharpoonup u \;\;\;\text{ weakly star in } 
  L^{\infty}(0,T;V)\text{ and weakly in }H^1(0,T;H),\label{conv2}\\
     &u_t^\eps\rightharpoonup u_t \;\;\;\text{ weakly star in } BV(0,T;X),\label{conv3}\\
       &\beta^\eps(u^\eps)\rightharpoonup\eta\;\;\;\text{ weakly in } \mathcal V'.\label{conv3a}
\end{align}
Here and below all convergence relations are implicitly intended to hold up
to extraction of a (non relabeled) subsequence of $\eps\searrow 0$.

Thanks to \eqref{conv1}-\eqref{conv2} and \eqref{est1*} we also infer
\begin{align}\label{conv4}
 &u^\eps(t)\rightharpoonup u(t) \;\;\;\text{ weakly in }V\;\;\;\text{ for all }t\in[0,T].
\end{align}
Next, condition \eqref{conv1} implies, thanks to the Aubin-Lions lemma, that
\begin{align}\label{conv4*}
  & u^\eps\rightarrow u \;\;\;\text{ strongly in } L^2(0,T;V).
\end{align}
A generalized version of the same lemma \cite[Cor.~4, Sec.~8]{Si}
implies, thanks to \eqref{conv2} and \eqref{conv3},
\begin{align}\label{conv5*}
  & u_t^\eps\rightarrow u_t \;\;\;\text{ strongly in } L^2(0,T;V').
\end{align}
From \eqref{conv3} and a proper version of the 
Helly selection principle \cite[Lemma 7.2]{DMDSM}, we infer
\begin{align}
  u_t^\eps(t)\rightharpoonup u_t(t) \;\;\;\text{ weakly in } X\;\;\;\;
   \text{ for all }t\in[0,T].\label{conv6} 
\end{align}
%
%
Combining this with \eqref{est2c}, we obtain more precisely
\begin{align}
 u_t^\eps(t)\rightarrow u_t(t) \;\;\;\text{ weakly in }V'
  \;\text{ and strongly in }D(A^{-1})\;\;\;\;\text{ for all }t\in[0,T].\label{conv7} 
\end{align}
\end{subequations}
Hence, using \eqref{conv}, we can take $\eps\searrow0$ in~\eqref{eq:wee} 
and get back~\eqref{eq:w}. Indeed, it is not difficult to check that all terms
pass to the limit. Notice however that, 
in view of \eqref{conv3a}, the {\sl $L^2$-scalar product}\/
$\dual{\beta^\eps(u^\eps),\varphi}$ is replaced by 
the $\calV'$-$\calV$ duality $\duav{\eta,\varphi}$ in the limit.

Let us now consider the weak formulation on subintervals. Taking $\varphi\in\mathcal V_t$, 
$t\in [0,T]$, we may rearrange terms in \eqref{eq:weet} to get
\begin{align}\nonumber
  & \dual{\beta^\eps(u^\eps),\varphi}_t
   = - \alpha(u^\eps_{t}(t),\varphi(t))_*
   + \alpha (u^\eps_1,\varphi(0))_*
   + \alpha\dual{u^\eps_{t},\varphi_t}_{*,t}\\
 \label{eq:weet2}
  & \mbox{}~~~~~
   - \dual{u^\eps_t,\varphi}_{*,t}
   - \delta\dual{u^\eps_t,\varphi}_t
   - \dual{A^{1/2} u^\eps, A^{1/2}\varphi}_t
   + \lambda\dual{u^\eps,\varphi}_t=0.     
\end{align}
Now, {\sl without extracting further subsequences}, it can be checked that,
as a consequence of \eqref{conv}, the right hand side tends to 
\begin{align}\nonumber
  & - \alpha(u_t(t),\varphi(t))_*
   + \alpha (u_1,\varphi(0))_*
   + \alpha \dual{u_t,\varphi_t}_{*,t}
   - \dual{u_t,\varphi}_{*,t}\\
 \label{eq:weet3}
  & \mbox{}~~~~~
   - \delta\dual{u_t,\varphi}_t
   - \dual{A^{1/2} u, A^{1/2}\varphi}_t
   + \lambda\dual{u,\varphi}_t
   =: \duav{\eta_{(t)},\varphi}_t.
\end{align}
%
Hence we have proved \eqref{eq:w} and \eqref{eq:wt}. The compatibility
property~\eqref{compa} is also a straighforward consequence
of this argument.

Next, to prove \eqref{graph}, according to~\eqref{minty}, we need to show
\begin{equation} \label{semi:1}
  \limsup_{\eps\searrow0} \,\duav{\beta^\eps(u^\eps),u^\eps}\leq\duav{\eta,u}.
\end{equation}
Thanks to \eqref{eq:wee} with $\varphi=u^\eps$, we have
\begin{align}
  & \duav{\beta^\eps(u^\eps),u^\eps} 
   = - \alpha (u^\eps_{t}(T),u^\eps(T))_*
    + \alpha (u^\eps_1,u^\eps_0)_*
    + \alpha \dual{u^\eps_{t},u^\eps_t}_*\nonumber\\
 \label{semi:2}
  & \mbox{}~~~~~ - \dual{u^\eps_t,u^\eps}_*
    - \delta \dual{u^\eps_t,u^\eps}
    - \|A^{1/2} u^\eps\|_{L^2(0,T;H)}^2
    + \lambda\|u^\eps\|_{L^2(0,T;H)}^2.
\end{align}
Then, we take the $\limsup$ of the above expression 
as $\eps\searrow0$. Then, using relations \eqref{conv}
and standard lower semicontinuity arguments
we infer that the $\limsup$ of the above expression
is less or equal than
\begin{align}
  & - \alpha (u_{t}(T),u(T))_*
    + \alpha (u_1,u_0)_*
    + \alpha \dual{u_{t},u_t}_*\nonumber
    - \dual{u_t,u}_*\\
 \label{semi:3}
  & \mbox{}~~~~~  - \delta \dual{u_t,u}
    - \|A^{1/2} u\|_{L^2(0,T;H)}^2
    + \lambda\|u\|_{L^2(0,T;H)}^2
    = \duav{\eta,u}, 
\end{align}
the last equality following from \eqref{eq:w} with the choice
$\varphi=u$. Combining \eqref{semi:2} with \eqref{semi:3} we obtain
\eqref{semi:1}, whence the first of \eqref{graph}. 
The same argument applied to the subinterval $(0,t)$ entails
$\eta_{(t)}\in \beta_w(u\llcorner_{(0,t)})$, for all $t\in(0,T]$, as desired.

\smallskip

Finally, we need to prove the energy {\sl inequality}\/ \eqref{energ:eq}.
To this aim, we consider the approximate energy balance~\eqref{energ:eps}
and take its $\liminf$ as $\eps\searrow0$.

Then, by standard lower semicontinuity arguments, it is clear that the
left hand side of \eqref{energ:eq} is less or equal than the 
$\liminf$ of the left hand side of \eqref{energ:eps}.
The more delicate point stands, of course, in dealing with 
the right hand sides. Indeed, we claim that there exists the limit
\begin{align}\no
  & \lim_{\eps\searrow 0} \Big( \frac{\alpha}{2}\|u_t^\eps(t_1)\|_{V'}^2
    + \frac{1}{2}\|A^{1/2} u^\eps(t_1)\|^2_{H}
    + J^\eps(u^\eps(t_1)) 
    - \frac{\lambda}{2}\|u^\eps(t_1)\|^2_H \Big) \\
 \label{semi:4}
  & \mbox{}~~~~~  
   = \Big( \frac{\alpha}{2}\|u_t(t_1)\|_{V'}^2
    + \frac{1}{2}\|A^{1/2} u(t_1)\|^2_{H}
    + J(u(t_1)) 
    - \frac{\lambda}{2}\|u(t_1)\|^2_H \Big),
\end{align}
at least for almost every $t_1 \in [0,t)$, surely including $t_1=0$.
We just sketch the proof of this fact, which follows closely the lines of 
the argument given in \cite[Section 3]{BRSS}
to which we refer the reader for more details.

First, we observe that the last summand passes to the limit
in view of \eqref{conv4} and the compact embedding $V\subset H$.
Next, the convergence
$$
 \Big( \frac{\alpha}{2}\|u_t^\eps(t_1)\|_{V'}^2
    + \frac{1}{2}\|A^{1/2} u^\eps(t_1)\|^2_{H} \Big) 
    \to \Big( \frac{\alpha}{2}\|u_t(t_1)\|_{V'}^2
    + \frac{1}{2}\|A^{1/2} u(t_1)\|^2_{H} \Big)
$$
holds for {\sl almost}\/ every choice of $t_1$ and up to extraction
of a further subsequence of $\eps\searrow 0$ in view of 
\eqref{conv4*} and \eqref{conv5*} (indeed, because these are just
$L^2$-bounds with respect to time, we cannot hope to get convergence
for {\sl every}\/ $t_1\in[0,T)$). Finally, we need to show
$$
  J^\eps(u^\eps(t_1))\rightarrow J(u(t_1)).
$$
This is the most delicate part, which proceeds exactly as 
in~\cite[Section 3]{BRSS}, to which 
the reader is referred. Note, finally, that \eqref{semi:4}
for $t_1=0$ can be easily proved as a direct consequence
of Lemma~\ref{iniz_data} (again, we refer the reader 
to~\cite{BRSS} for details). The proof is concluded.


\end{document}